\documentclass[12pt]{article}
\usepackage[utf8]{inputenc}
\usepackage{graphicx}
\usepackage[hidelinks]{hyperref}
\hypersetup{
	colorlinks=true,
	linkcolor=blue,
	citecolor=black,
	filecolor=magenta,
	urlcolor=blue, }

\usepackage{cite, url}
\usepackage{amsfonts,float,lscape,color,amsmath,amsthm,amssymb,mathrsfs,mathtools,bbm, faktor}
\usepackage{cleveref}
\usepackage[dvipsnames, table]{xcolor}

\usepackage{pgfplots,pgf}
\pgfplotsset{compat=newest}
\usetikzlibrary{intersections,hobby,patterns, patterns.meta,shadings,arrows,decorations,backgrounds}
\pgfplotsset{my style/.append style={axis x line=middle, axis y line=
		middle, xlabel={$x$}, ylabel={$y$}, axis equal }}
\usepackage{tikz,tikz-cd,tkz-fct}

\newcommand{\R}{\mathbb{R}} 
\newcommand{\C}{\mathbb{C}}
\newcommand{\T}{\mathbb{T}}
\newcommand{\Z}{\mathbb{Z}}

\theoremstyle{plain}
\newtheorem{theorem}{Theorem}[section]
\newtheorem*{theorem*}{Theorem}
\newtheorem{prop}[theorem]{Proposition}
\newtheorem{lemma}[theorem]{Lemma}
\newtheorem{cor}[theorem]{Corollary}
\theoremstyle{definition} 
\newtheorem{defn}[theorem]{Definition}
\newtheorem{ex}[theorem]{Example}
\theoremstyle{remark}
\newtheorem{rem}[theorem]{Remark}

\theoremstyle{plain}
\newtheorem*{namedthm}{\namedthmname}
\newcounter{namedthm}



\usepackage{upgreek}
\usepackage{tabularray,longtable}

\title{Subspaces of $L^2(\R^n)$ Invariant Under Shifts by a Crystal Group}
\author{Tom Potter and Keith Taylor}
\date{\today}

\makeatletter
\def\namedlabel#1#2{\begingroup
	#2%
	\def\@currentlabel{#2}%
	\phantomsection\label{#1}\endgroup
}
\makeatother

\begin{document}
	
\maketitle
		\begin{abstract}
For a crystal group $\Gamma$ in dimension $n$, a closed subspace $\mathcal{V}$ of $L^2(\R^n)$ is called 
$\Gamma$--shift-invariant if, for every $f\in\mathcal{V}$, the shifts of $f$ by every element of $\Gamma$
also belong to $\mathcal{V}$. The main purpose of this paper is to provide a characterization of
the $\Gamma$--shift-invariant closed subspaces of $L^2(\R^n)$ that is valid for all crystal groups, 
non-symmorphic as well as symmorphic.
	\end{abstract}

	\section{Introduction} 

Let $\Gamma$ be a crystal group in dimension $n$ with point group $\Pi$ and associated lattice
$\mathcal{L}$. 
Detailed definitions are given in Section 3. Elements of $\Gamma$ are written $[x,M]$, where
$M$ is an orthogonal matrix that belongs to $\Pi$ and $x$ is a vector in $\R^n$.  The crystal group is called 
{\em symmorphic} if $[0,M]\in\Gamma$, for all $M\in\Pi$; or, equivalently, if $\Gamma$ is isomorphic to the
semi-direct product $\mathcal{L}\rtimes\Pi$. If this is not the case, then $\Gamma$ is {\em non-symmorphic}.

For a function $f$
on $\R^n$, the {\em shift} of $f$ by $[x,M]$ is $\pi[x,M]f$ given by
\[
\pi[x,M]f(z)=f\left(M^{-1}z-x\right),
\]
for $z\in\R^n$.
When these shifts are applied to functions in the Hilbert space of all square-integrable functions, $L^2(\R^n)$, 
they constitute a unitary
representation $\pi$ of $\Gamma$, which we call the natural representation. 
A closed subspace $\mathcal{V}$ of $L^2(\R^n)$ is called $\pi$-invariant,
or {\em$\Gamma$--shift-invariant}, if $\pi[x,M]f\in\mathcal{V}$, for every $f\in\mathcal{V}$ and every
$[x,M]\in\Gamma$. Our purpose is to provide a characterization of the $\Gamma$--shift-invariant closed
subspaces of $L^2(\R^n)$.

The term shift-invariant closed subspace means a closed subspace $\mathcal{V}$ of $L^2(\R^n)$ such that 
$T_kf\in\mathcal{V}$, for all $f\in\mathcal{V}$ and all $k\in\Z^n$, where 
$T_kf(z)=f(z-k)$, for $z\in\R^n$. In every dimension, there is one
crystal group that is abelian and it is the set of shifts by lattice points; so shift-invariance is
$\Gamma$--shift invariance when $\Gamma$ is isomorphic to $\Z^n$.
Understanding shift-invariant subspaces became
more important with the rise of wavelet analysis (see \cite{Daubechies} for a general introduction to 
the classical theory of wavelets) as the central subspace in a multiresolution analysis is shift-invariant; 
see also \cite{RonShen} for more on the role of shift-invariant subspaces in wavelet theory.
A characterization of shift-invariant subspaces was known in 1964 by Helson \cite{Helson} and was
 reformulated in 2000 by Bownik \cite{Bownik} where he applied it to wavelet theory. 

The theory of wavelets has many generalizations and variations. One generalization, wavelets with
composite dilations, was introduced in \cite{Guo1} and \cite{Guo2}. The term \textit{composite dilations} refers
to the addition of another set $\mathcal{B}$ of matrices that are moving functions, besides the shifts by
vectors in a lattice and the powers of a single dilating matrix. Often, $\mathcal{B}$ is a finite group
of measure preserving matrices. Blanchard developed the theory in this case in \cite{Blanchard1}; see
also \cite{Blanchard2} where the focus was on developing Haar-type wavelet systems when $\mathcal{B}$
leaves the shift lattice invariant and the semi-direct product of $\mathcal{B}$ with the lattice is actually
a crystal group, necessarily symmorphic.
 In \cite{KeithJosh} MacArthur and one of the current authors introduced the concept of
a multiresolution analysis where the shifts came from an arbitrary crystal group and cast the theory in
the context of the abstract approach of \cite{Baggett2}. Independently, Gonz\'{a}lez and Moure \cite{Gonzalez}
also formulated the theory for shifts by a crystal group. In the definition of a multiresolution analysis or a
generalized multiresolution analysis for either composite dilations, with a finite $\mathcal{B}$,
 or shifts by a crystal group, the central subspace in the multiresolution analysis is left invariant under
shifts by a crystal group $\Gamma$.

Another reason to be concerned about the nature of $\Gamma$--shift-invariant closed subspaces is 
the emergence of topological quantum chemistry  (see \cite{Cano2} and \cite{Cano1}, for example) where
band representations play a significant role. A band representation of a crystal group $\Gamma$ is the
restriction of the natural representation $\pi$ to a  closed $\Gamma$--shift-invariant subspace generated
by a single function in $L^2(\R^n)$. A readable review of band representations can be found in \cite{Bacry}.

In 2020, Barbieri, et al, \cite{Barbieri} considered $\Gamma$--shift-invariance of
closed subspaces of $L^2(\mathcal{G})$, where $\mathcal{G}$ is a locally compact abelian group and
$\Gamma$ is the semi-direct product, $\Lambda\rtimes G$, where $\Lambda$ is a uniform lattice in
$\mathcal{G}$ and $G$ is a finite group of continuous automorphisms of $\mathcal{G}$ that fix
$\Lambda$. Theorem 4.7 of \cite{Barbieri} provides a characterization of $\Gamma$--shift-invariant closed 
subspaces of $L^2(\mathcal{G})$. When $\mathcal{G}=\R^n$, this generalizes the characterization of 
shift-invariant closed subspaces of $L^2(\R^n)$ given in \cite{Bownik}  to $\Gamma$--shift-invariance,
where $\Gamma$ is a symmorphic crystal group. 

Our main result, given in Theorem \ref{main_theorem} and
reformulated in Corollary \ref{main_cor}, applies to all crystal groups. Thus, the contribution of this paper
is to describe the situation for non-symmorphic crystal groups. The two crystal groups in one dimension are
symmorphic, while 4 of the 17 two dimensional crystal groups are non-symmorphic. In three dimensions, 
146 out of the 219 different (up to isomorphism) crystal groups are non-symmorphic. It gets more dramatic 
as the dimension increases. In dimensions 5 and 6, more than 85\% and 97\%, respectively, of the crystal
groups are non-symmorphic (see \cite{Plesken} for exact numbers). A readable account of the process of
enumeration of crystal groups is given in \cite{Opg}. Also, Section 2.5 of \cite{Opg} contains an explanation of 
the connection between a vector system for a given point group and the non-symmorphic crystal groups
with that point group. We will discuss how a vector system can be defined for a particular $\Gamma$ in
Section 3 and it turns out that the vector system plays a key role in our main result.

Imprecisely stated, we show, in
Theorem \ref{main_theorem} and Corollary \ref{main_cor}, that
there is a one-to-one correspondence between $\Gamma$--shift-invariant closed subspaces of
$L^2(\R^n)$ and certain maps (twisted range functions, where the twist comes from the vector system) 
from an open subset of $\R^n$ and the
closed subspaces of $\ell^2(\mathcal{L}^*)$, where $\mathcal{L}^*$ is the dual lattice of $\mathcal{L}$.
 If $\Gamma=\Z^n$, then
the characterization in Theorem \ref{main_theorem} reduces to the characterization of shift-invariant
subspaces given in Proposition 1.5 in \cite{Bownik}. When $\Gamma$ happens to be symmorphic,
Theorem \ref{main_theorem} is equivalent to the characterization for symmorphic crystal groups found
in \cite{Barbieri}. Our approach draws on the presentation of 
range functions by Bownik and Ross in \cite{BownikRoss}, to the extent that $\Gamma$--shift-invariance
implies invariance under shifts by the lattice $\mathcal{L}$. We view $\pi$ as a unitary representation of
$\Gamma$ and our strategy is to use unitary maps to transform $\pi$ into an equivalent representation
that can be analyzed using range function methods.

We gather together preliminary facts and known results in Section 2 while Section 3 is devoted to organizing
the properties of crystal groups that are needed. The unitary maps that will be used to transform $\pi$
are introduced in Section 4 and applied in the next section. In the final section, we prove our main theorem,
Theorem \ref{main_theorem}. We also present an example of a non-symmorphic wallpaper group
$\Gamma$ and a $\Gamma$--shift-invariant closed subspace of $L^2(\R^2)$ where the impact of the
non-symmorphic nature of this $\Gamma$ is illustrated and the twist coming from the vector system is clear.

\medskip\noindent
{\bf Acknowledgements}: Our interest in finding a characterization of $\Gamma$--shift-invariant closed 
subspaces originated with a question posed by Azita Mayeli. We also thank a referee of an earlier version of this
paper who brought \cite{Barbieri} to our attention.

\section{Preliminaries}
Let $n$ be a positive integer and let $\R^n$ denote Euclidean space with $x\in\R^n$ being considered as
a column vector. For $1\leq p\leq\infty$, $L^p(\R^n)$ denotes the standard Lebesgue space with respect to
Lebesgue measure on $\R^n$ and $\|\cdot\|_p$ denotes the usual norm in $L^p(\R^n)$. We use the following version of the Fourier transform. For $f\in L^1(\R^n)$ and $\omega\in\R^n$, let
$\widehat{f}(\omega)=\int_{\R^n}f(x)e^{2\pi i \omega\cdot x}dx$. 
 By Plancherel's Theorem, there is a unitary map $\mathcal{F}$ on $L^2(\R^n)$
such that $\mathcal{F}f=\widehat{f}$, for all $f\in L^1(\R^n)\cap L^2(\R^n)$. We call both $\mathcal{F}$ and 
the map $f\mapsto\widehat{f}$ the Fourier transform. We will often use $\widehat{f}$ instead of $\mathcal{F}f$, for
any $f\in L^2(\R^n)$.

If we consider $\R^n$ as a topological group with addition of vectors as the group operation, a subgroup 
$\mathcal{L}$ of $\R^n$ is called a lattice in $\R^n$ if $\mathcal{L}$ is discrete and $\R^n/\mathcal{L}$ is compact.
The dual lattice of $\mathcal{L}$ is $\mathcal{L}^*=\{\nu\in\R^n:\nu\cdot k\in\Z\}$. That is, $\nu\in\R^n$ is 
in $\mathcal{L}^*$ if and only if $e^{2\pi i\nu\cdot k}=1$, for all $k\in\mathcal{L}$. 
If $\mathcal{L}$ is a lattice in $\R^n$, then there exists an invertible $n\times n$ matrix $B$ such that
$\mathcal{L}=B\Z^n=\{Bk:k\in\Z^n\}$. Then $\mathcal{L}^*=(B^{-1})^t\Z^n$, where $(B^{-1})^t$ is the
transpose of the inverse of $B$. Let $Q=\{(B^{-1})^t\theta:\theta\in\left[-\frac{1}{2},\frac{1}{2}\right)^n\}$. 
Then $\R^n$ is the 
disjoint union of the $Q+\nu$, $\nu\in\mathcal{L}^*$. We equip $Q$ with the restriction of Lebesgue 
measure.

We will need a form of the Poisson summation formula. For $g\in L^1(\R^n)$, the periodization
$\mathcal{P}_{\mathcal{L}^*}g$ of $g$ with respect to $\mathcal{L}^*$ is given by 
\[
\mathcal{P}_{\mathcal{L}^*}g(\theta)=
\sum_{\nu\in\mathcal{L}^*}g(\theta+\nu),
\]
for any $\theta\in\R^n$ for which the series converges.
\begin{lemma}\label{Poisson1}
Let $\mathcal{L}$ be a lattice in $\R^n$ and let $\mathcal{L}^*$ be its dual lattice. Let $g\in L^1(\R^n)$. 
Then $\sum_{\nu\in\mathcal{L}^*}g(\theta+\nu)$, converges for a.e. $\theta\in\R^n$ and 
$\mathcal{P}_{\mathcal{L}^*}g\big|_Q\in L^1(Q)$.
\end{lemma}
\begin{proof}
We adapt the standard argument as in the proof of Theorem 8.31 of \cite{FollandReal}.
\begin{eqnarray*}
\int_Q\sum_{\nu\in\mathcal{L}^*}|g(\theta+\nu)|d\theta & = &
\sum_{\nu\in\mathcal{L}^*}\int_Q|g(\theta+\nu)|d\theta\\
& = &
\sum_{\nu\in\mathcal{L}^*}\int_{Q+\nu}|g(\theta)|d\theta=\int_{\R^n}|g(\theta)|d\theta=\|g\|_1.
\end{eqnarray*}
By the sum variation of the Dominated Convergence Theorem, 
$\sum_{\nu\in\mathcal{L}^*}g(\theta+\nu)$ converges for a.e. $\theta\in Q$ and 
$\mathcal{P}_{\mathcal{L}^*}g\big|_Q\in L^1(Q)$. 
Since $\mathcal{P}_{\mathcal{L}^*}g$ is $\mathcal{L}^*$-periodic,
 $\sum_{\nu\in\mathcal{L}^*}g(\theta+\nu)$ converges for a.e. $\theta\in \R^n$.
\end{proof}
\begin{prop}\label{Poisson}
Let $\mathcal{L}$ be a lattice in $\R^n$ and let $\mathcal{L}^*$ be its dual lattice. Let $g\in L^1(\R^n)$. 
For each $\ell\in\mathcal{L}$, 
\[
\int_Q \mathcal{P}_{\mathcal{L}^*}g(\theta)e^{2\pi i\theta\cdot\ell}d\theta=\widehat{g}(\ell).
\]
Moreover, if $\widehat{g}(\ell)=0$, for all $\ell\in\mathcal{L}$, then $\mathcal{P}_{\mathcal{L}^*}g=0$.
\end{prop}
\begin{proof}
Again, dominated convergence justifies the interchange of sum and integral below. For $\ell\in\mathcal{L}$,
\begin{eqnarray*}
\int_Q \mathcal{P}_{\mathcal{L}^*}g(\theta)e^{2\pi i\theta\cdot\ell}d\theta &=&
\int_Q \sum_{\nu\in\mathcal{L}^*}g(\theta+\nu)e^{2\pi i\theta\cdot\ell}d\theta=
 \sum_{\nu\in\mathcal{L}^*}\int_Qg(\theta+\nu)e^{2\pi i\theta\cdot\ell}d\theta\\
&=&
 \sum_{\nu\in\mathcal{L}^*}\int_{Q+\nu}g(\theta)e^{2\pi i(\theta-\nu)\cdot\ell}d\theta=
\sum_{\nu\in\mathcal{L}^*}\int_{Q+\nu}g(\theta)e^{2\pi i\theta\cdot\ell}d\theta\\
&=&
\int_{\R^n}g(\theta)e^{2\pi i\theta\cdot\ell}d\theta=\widehat{g}(\ell).
\end{eqnarray*}
For the last statement in the proposition, we do a linear change of variables. Let $C=(B^{-1})^t$, where
$\mathcal{L}=B\Z^n$. That is,
$C$ is the invertible matrix such that $\mathcal{L}^*=C\Z^n$.
Now define $h\in L^1(\R^n)$ by $h(\omega)=g(C\omega)$, for a.e. $\omega\in\R^n$. Form the $\Z^n$
periodization of $h$. That is,
\[
\mathcal{P}_{\Z^n}h(\omega)=\sum_{k\in\Z^n}h(\omega+k)=\sum_{k\in\Z^n}g(C\omega+Ck)=
\mathcal{P}_{\mathcal{L}^*}g(C\omega),
\]
which converges for a.e. $\omega\in\R^n$, and $\mathcal{P}_{\Z^n}h$ is integrable over 
$\left[-\frac{1}{2},\frac{1}{2}\right)^n$. Let $(c_j)_{j\in\Z^n}$ be the Fourier multi-series of the
periodic function $\mathcal{P}_{\Z^n}h$. For each $j\in\Z^n$, compute the Fourier coefficient $c_j$ using the
change of variables $\theta=C\omega$ and noting that $C^{-1}=B^t$. So
\begin{eqnarray*}
c_j &=& 
\int_{\left[-\frac{1}{2},\frac{1}{2}\right)^n}\mathcal{P}_{\Z^n}h(\omega)e^{2\pi i\omega\cdot j}d\omega=
\int_{\left[-\frac{1}{2},\frac{1}{2}\right)^n}\mathcal{P}_{\mathcal{L}^*}g(C\omega)e^{2\pi i\omega\cdot j}d\omega\\
&=&
|\det(B)|\int_Q\mathcal{P}_{\mathcal{L}^*}g(\theta)e^{2\pi i\theta\cdot Bj}d\theta=|\det(B)|\widehat{g}(Bj).
\end{eqnarray*}
Now $Bj\in\mathcal{L}^*$, for all $j\in\Z^n$. So, if $\widehat{g}(\ell)=0$, for all $\ell\in\mathcal{L}^*$, then
every Fourier coefficient of $\mathcal{P}_{\Z^n}h$ is 0. This implies $\mathcal{P}_{\Z^n}h=0$, which implies 
$\mathcal{P}_{\mathcal{L}^*}g=0$.
\end{proof}

\begin{cor}\label{determining}
With the notation of Proposition \ref{Poisson}, let $Y$ be an open subset of $\R^n$ such that
$(Y+\nu)\cap(Y+\kappa)=\emptyset$ if $\nu,\kappa\in\mathcal{L}^*$ with $\nu\neq\kappa$.
Let $g\in L^1(\R^n)$ be such that $g(\theta)=0$, for a.e. $\theta\in\R^n\setminus Y$. If
$\widehat{g}(\ell)=0$, for all $\ell\in\mathcal{L}$, then $g=0$.
\end{cor}
\begin{proof}
Under the assumption on $g$, $\mathcal{P}_{\mathcal{L}^*}g$ agrees with $g$ on $Y$, so 
Proposition \ref{Poisson} implies $g=0$.
\end{proof}
\begin{rem}\label{determing2}
For each $\ell\in\mathcal{L}$, let $e_\ell(\theta)=e^{2\pi i\theta\cdot\ell}$, for all $\theta\in\R^n$. In the
language of \cite{BownikRoss}, Definition 2.2, $\mathcal{D}=\{e_\ell|_{_Y}:\ell\in\mathcal{L}\}$ is
a determining set for $L^1(Y)$.
\end{rem}

The concept of a range function goes back to \cite{Helson} and was used in \cite{Bownik} for the characterization
of closed subspaces of $L^2(\R^n)$ that are shift-invariant. The treatment of range functions given in
\cite{BownikRoss} is most useful for our purposes.
 For a Hilbert space $\mathcal{H}$, let $Gr(\mathcal{H})$ denote the set of all closed
subspaces of $\mathcal{H}$, the Grassmannian of $\mathcal{H}$. If $\mathcal{K}\in Gr(\mathcal{H})$, let $P_{\mathcal{K}}$ denote the
orthogonal projection of $\mathcal{H}$ onto $\mathcal{K}$.
\begin{defn}
Let $\mathcal{H}$ be a separable Hilbert space and let 
$(\mathcal{X},\Sigma)$ be a measurable space. A measurable {\em range function} for $\mathcal{H}$
based on $\mathcal{X}$, is a mapping $J:\mathcal{X}\to Gr(\mathcal{H})$ such that, for any $\xi,\eta\in
\mathcal{H}$, the map $x\to\langle P_{J(x)}\xi,\eta\rangle$ is measurable.
\end{defn}
In all the cases we consider here, $\mathcal{H}$ is a separable Hilbert space. A map $F:\mathcal{X}\to
\mathcal{H}$ is measurable if $x\to\langle F(x),\eta\rangle$ is a measurable complex-valued function, for
each $\eta\in\mathcal{H}$.
If $\mu$ is a positive measure on $(\mathcal{X},\Sigma)$, let 
\[
L^2(\mathcal{X},\mu,\mathcal{H})=\left\{F:\mathcal{X}\to\mathcal{H}\big| F \text{ is measurable and }
\int_{\mathcal{X}}\|F(x)\|^2d\mu(x)<\infty\right\},
\]
with the usual identification of functions that agree $\mu$-almost everywhere. 
If the measure $\mu$ is
understood from context, we write $L^2(\mathcal{X},\mathcal{H})$.
The inner product of 
$F_1,F_2\in L^2(\mathcal{X},\mathcal{H})$ is given by
\[
\langle F_1,F_2\rangle=
\int_{\mathcal{X}}\langle F_1(x),F_2(x)\rangle d\mu(x).
\]
 Let $J$ be
a measurable range function for $\mathcal{H}$ based on $\mathcal{X}$. Then, for any 
$F\in L^2(\mathcal{X},\mathcal{H})$, $x\to\langle\mathcal{P}_{J(x)}F(x),\eta\rangle$ is measurable, 
for every $\eta\in\mathcal{H}$. 
Given $J$, define $\mathcal{M}_J=\{F\in L^2(\mathcal{X},\mathcal{H}): F(x)\in J(x), \text{ for }\mu\text{-a.e. }
x\in\mathcal{X}\}$. The next Proposition gathers together Propositions 2.1, 2.2, and 2.3 of \cite{BownikRoss}.
\begin{prop}\label{Range}
Let $J$and $K$  be measurable range functions for $\mathcal{H}$ based on $\mathcal{X}$. Then

{\rm (i)} $\mathcal{M}_J$ is a closed subspace of $L^2(\mathcal{X},\mathcal{H})$,

{\rm (ii)} if $\mathcal{P}_{\mathcal{M}_J}$ is the orthogonal projection of $L^2(\mathcal{X},\mathcal{H})$ onto $\mathcal{M}_J$, then,
for any $F\in L^2(\mathcal{X},\mathcal{H})$, $\left(\mathcal{P}_{\mathcal{M}_J}F\right)(x)=\mathcal{P}_{J(x)}F(x)$, for
a.e. $x\in\mathcal{X}$, and

{\rm (iii)} $\mathcal{M}_J=\mathcal{M}_K$ if and only if $J(x)=K(x)$, for a.e. $x\in \mathcal{X}$.
\end{prop}

If $\mathcal{H}$ is a Hilbert space, let $\mathcal{U}(\mathcal{H})$ denote the group of all unitary operators
on $\mathcal{H}$. If $\Lambda$ is any discrete group, a unitary representation of $\Lambda$ on $\mathcal{H}$
 is a
homomorphism $\sigma:\Lambda\to\mathcal{U}(\mathcal{H})$. Let $\sigma_1$ and $\sigma_2$ be unitary representations of $\Lambda$ on $\mathcal{H}_1$ and $\mathcal{H}_2$, respectively. If there exists
a unitary  map $W:\mathcal{H}_1\to\mathcal{H}_2$ such that $W\sigma_1(a)=\sigma_2(a)W$, for all
$a\in\Lambda$, then $\sigma_1$ and $\sigma_2$ are called equivalent.

\section{Crystal groups}

Let ${\rm O}_n$ denote the compact group of orthogonal $n\times n$ real matrices. 
For $x\in\R^n$ and $A\in{\rm O}_n$, define $[x,A]$ to be the affine map, $z\to A(z+x)$, of $\R^n$. Let
\[
{\rm Iso}_n(\R)=\{[x,A]:x\in\R^n, A\in {\rm O}_n\},
\]
the group of all isometries of $\R^n$. The composition of isometries is the group product: For
$[x,A],[y,B]\in{\rm Iso}_n(\R)$, $ [x,A][y,B]=[B^{-1}x+y,AB]$ and $[x,A]^{-1}=[-Ax,A^{-1}]$. Note that
$[0,{\rm id}]$ is the identity element of ${\rm Iso}_n(\R)$, where ${\rm id}$ denotes the identity $n\times n$ matrix.
Given the product topology of $\R^n\times{\rm O}_n$, ${\rm Iso}_n(\R)$ is a locally compact group.
Let ${\rm Trans}_n=\{[x,{\rm id}]:x\in\R^n\}$, the set of pure translations in ${\rm Iso}_n(\R)$. This is
a closed normal subgroup of ${\rm Iso}_n(\R)$. Let $q:{\rm Iso}_n(\R)\to{\rm O}_n$ be the homomorphism 
given by
$q[x,A]=A$, for all $[x,A]\in{\rm Iso}_n(\R)$.

A crystal group is a discrete subgroup $\Gamma$ of ${\rm Iso}_n(\R)$ such that $\R^n/\Gamma$ is
compact, where $\R^n/\Gamma$ is the set of all $\Gamma$-orbits, with the quotient topology. Section 7.5 of \cite{Ratcliffe} presents the basic properties of crystal groups. The translation
subgroup of $\Gamma$ is $T=\Gamma\cap {\rm Trans}_n$, which is a normal subgroup,
 and the point group is $\Pi=q(\Gamma)$, which is isomorphic to the quotient group $\Gamma/T$.
Then $\Pi$ is a finite subgroup of ${\rm O}_n$, while $T$ is a free  abelian group of rank $n$. Indeed,
$\mathcal{L}=\{\ell\in\R^n:[\ell,{\rm id}]\in T\}$ is a lattice in $\R^n$. For $[\ell,{\rm id}]\in T$ and $[x,M]\in\Gamma$,
\[
[x,M][\ell,{\rm id}][x,M]^{-1}=[x+\ell,M][-Mx,M^{-1}]=[M\ell,{\rm id}].
\]
This shows that if $\ell\in\mathcal{L}$, then $M\ell\in\mathcal{L}$, for all $M\in\Pi$.

Fix a cross-section $\gamma:\Pi\to\Gamma$ of the $T$-cosets in  $\Gamma$; so $q\big(\gamma(M)\big)
=M$, for all $M\in\Pi$. With $\gamma$ fixed, for each $M\in\Pi$, let $x_M\in\R^n$ be such that 
$\gamma(M)=[x_M,M]$. 
Then $\Gamma=\{[\ell+x_M,M]:\ell\in\mathcal{L},M\in\Pi\}$. For many crystal groups, $[0,M]\in
\Gamma$, for all $M\in\Pi$. Then, choose $x_M=0$, for each $M\in\Pi$. 
When this can be done, $\Gamma$ is isomorphic to
the semidirect product $\mathcal{L}\rtimes\Pi$. As pointed out in the introduction, such
 crystal groups are called symmorphic.

\begin{defn}
Let $\Gamma$ be a crystal group in dimension $n$ with point group $\Pi$. 
A {\em vector system associated with $\Gamma$} is a map $M\mapsto x_M$
of $\Pi$ into $\R^n$ such that $[x_M,M]\in\Gamma$, for all $M\in\Pi$.
\end{defn}

\begin{lemma}\label{cocycle} If $\{x_M:M\in\Pi\}$ is a vector system associated with $\Gamma$, then
$N^{-1}x_K+x_N-x_{KN}\in\mathcal{L}$, for all $K,N\in\Pi$.
\end{lemma}
\begin{proof}
For $K,N\in\Pi$, calculate the following element of $\Gamma$:
\[
\begin{array}{ll}
[x_K,K][x_N,N][x_{KN},KN]^{-1} & =[N^{-1}x_K+x_N,KN][-KNx_{KN},(KN)^{-1}]\\
& = [Kx_K+KNx_N-KNx_{KN},{\rm id}]\in T.
\end{array}
\]
Thus, $Kx_K+KNx_N-KNx_{KN}\in\mathcal{L}$. Then
\[
N^{-1}x_K+x_N-x_{KN}=N^{-1}K^{-1}\big(Kx_K+KNx_N-KNx_{KN}\big)
\]
is in $\mathcal{L}$ as well.
\end{proof}

Associated with the lattice $\mathcal{L}$ is the dual lattice
$\mathcal{L}^*=\{\nu\in\R^n: \nu\cdot\ell\in\Z, \text{ for all }\ell\in\mathcal{L}\}$. 
 Since $M^t=M^{-1}$, for any $M\in\Pi\subseteq {\rm O}_n$, this implies
$\mathcal{L}^*$ is left invariant when multiplied by members of $\Pi$. With this action of
$\Pi$ on $\mathcal{L}^*$, form the semidirect product $\mathcal{L}^*\rtimes\Pi=\{(\nu,M):
\nu\in\mathcal{L}^*, M\in\Pi\}$, with the group product given by
\[
(\kappa,L)(\nu,M)=(M^{-1}\kappa+\nu,LM), \text{ for } (\kappa,L), (\nu,M)\in\mathcal{L}^*\rtimes\Pi.
\]
This auxiliary group might be isomorphic to $\Gamma$, but this is not always the case. As with $\Gamma$,
 $\mathcal{L}^*\rtimes\Pi$ acts on $\R^n$.  For $(\nu,M)\in \mathcal{L}^*\rtimes\Pi$ and $\omega\in \R^n$,
let $(\nu,M)\cdot\omega=M(\omega+\nu)$. This identifies $\mathcal{L}^*\rtimes\Pi$ with a discrete group
of isometries of $\R^n$ such that, since $\mathcal{L}^*$ is a full-rank lattice, $\R^n/(\mathcal{L}^*\rtimes\Pi)$ is
compact. That is, $\mathcal{L}^*\rtimes\Pi$ is also a crystal group. Let $\Gamma^*=\mathcal{L}^*\rtimes\Pi$.

For any $\omega\in\R^n$, let $\Gamma^*\omega=\{M(\omega+\nu):(\nu,M)\in\Gamma^*\}$, the 
$\Gamma^*$-orbit of $\omega$, and
let $\Gamma^*_\omega=\{(\nu,M)\in\Gamma^*:M(\omega+\nu)=\omega\}$, the stabilizer of $\omega$.
There exist points $\omega$ in $\R^n$ such that $\Gamma^*_\omega=\{(0,{\rm id})\}$ (see Theorem 6.6.12 of
\cite{Ratcliffe}).
Fix $\omega_0\in\R^n$ such that $\Gamma^*_{\omega_0}=\{(0,{\rm id})\}$. For each $(\nu,M)\in\Gamma^*
\setminus\{(0,{\rm id})\}$, let 
$H_{(\nu,M)}=\{\omega\in\R^n:\|\omega-\omega_0\|<\|\omega-M(\omega_0+\nu)\|\}$.
\begin{defn}\label{Dirichlet}
The {\em Dirichlet domain} for $\Gamma^*$ containing $\omega_0$ is
\[
\Omega_{\omega_0}=\bigcap\{H_{(\nu,M)}:(\nu,M)\in\Gamma^*, (\nu,M)\neq (0,{\rm id})\}.
\]
\end{defn}

Let $A\subseteq\R^n$, then  $A+\nu=\{\omega+\nu:\omega\in A\}$, for $\nu\in\mathcal{L}^*$, and
$MA=\{M\omega:\omega\in A\}$, for $M\in\Pi$. Also,
$\partial A$ denotes the boundary of $A$.

\begin{prop}\label{Rproperties}
Let $\omega_0\in\R^n$ be such that $\Gamma^*_{\omega_0}=\{(0,{\rm id})\}$ and let 
$\Omega=\Omega_{\omega_0}$.
Then $\Omega$ has the following properties:
\begin{enumerate}
\item $\Omega$  is open
\item $\Omega$  is convex
\item For $(\kappa,L),(\nu,M)\in\Gamma^*, (\kappa,L)\neq(\nu,M)$,
$ (L(\Omega+\kappa))\cap(M(\Omega+\nu))=\emptyset$
\item  $\cup_{(\nu,M)\in\Gamma^*}M(\overline{\Omega}+\nu)=\R^n$
\item $\cup_{(\nu,M)\in\Gamma^*}M(\Omega+\nu)$ is a co-null subset of $\R^n$.
\end{enumerate}
\end{prop}
\begin{proof}
Properties 1, 2, 3, and 4 are well-known (see, for example, Theorem 6.6.13 and the definition
of {\em fundamental domain} on page 233 of \cite{Ratcliffe}). Since $\Omega$ is open and convex
$\partial\Omega$ has Lebesgue measure 0. Each map $\omega\to M(\omega+\nu)$ is an isometry, so
$\partial(M(\Omega+\nu))=M(\partial\Omega+\nu)$ is a null set, for each $(\nu,M)\in\Gamma^*$.
Now $\R^n\setminus \cup_{(\nu,M)\in\Gamma^*}M(\Omega+\nu)\subseteq 
\cup_{(\nu,M)\in\Gamma^*}M(\partial\Omega+\nu)$, which is a null set since $\Gamma^*$ is countable.
This implies 5.
\end{proof}
Let $\Pi\Omega=\cup_{M\in\Pi}M\Omega$. Then $\Pi\Omega$ is a fundamental domain for $\mathcal{L}^*$.
In particular, we have the following corollary.
\begin{cor}\label{disjoint_translates}
If $\nu_1,\nu_2\in\mathcal{L}^*$, $\nu_1\neq\nu_2$, then $(\Pi\Omega+\nu_1)\cap(\Pi\Omega+\nu_2)=\emptyset$
and $\cup_{\nu\in\mathcal{L}^*}(\Pi\Omega+\nu)$ is co-null in $\R^n$.
\end{cor}

\section{Several Unitary Transformations}
Throughout this section, $\Gamma$ is a fixed crystal group with $\mathcal{L}$, $\mathcal{L}^*$, $\Pi$, $T$, 
$\Omega$, and a vector system $\{x_M:M\in\Pi\}$ as in Section 3.
In order to illuminate aspects of the natural unitary representation of the crystal group $\Gamma$
in the next section, the Hilbert space
$L^2(\R^n)$ will be transformed into a number of different realizations by specific unitary maps. The first is the
Fourier transform viewed as a unitary map of $L^2(\R^n)$ onto itself. 

For the other Hilbert spaces that arise, we equip the open subsets
$\Omega$ and $\Pi\Omega$ with the restriction of Lebesgue measure and $\Omega\times\Pi\times\mathcal{L}^*$
with the product measure, where $\Pi$ and $\mathcal{L}^*$ have counting measure.
We can use $\Omega\times\Pi\times\mathcal{L}^*$ to almost parametrize $\R^n$. Let 
$X=\cup_{(L,\nu)\in\Pi\times\mathcal{L}^*}L(\Omega+\nu)$. By Proposition \ref{Rproperties}, $X$ is an
open co-null subset of $\R^n$.

\begin{lemma}\label{phihom}
The map $\phi:\Omega\times\Pi\times\mathcal{L}^*\to X$ given by $\phi(\omega,L,\nu)=L(\omega+\nu)$ is a measure preserving homeomorphism.
\end{lemma}
\begin{proof}
Routine.
\end{proof}

For any $\theta\in X$, let $\omega_\theta\in\Omega$, $L_\theta\in\Pi$, and $\nu_\theta\in\mathcal{L}^*$ be such
that $L_\theta(\omega_\theta+\nu_\theta)=\theta$.
For each $\xi\in L^2(\R^n)$, let $W_1\xi=\xi\circ\phi$. Then $W_1\xi$ is a  measurable function on $\Omega\times\Pi\times\mathcal{L}^*$.

\begin{lemma}\label{W1}
The map
$W_1$ is a unitary from $L^2(\R^n)$ onto $L^2(\Omega\times\Pi\times\mathcal{L}^*)$ and 
$W_1^{-1}$ is given by, for 
$g\in L^2(\Omega\times\Pi\times\mathcal{L}^*)$, 
$W_1^{-1}g(\theta)=
g(\omega_\theta,L_\theta,\nu_\theta)$,  for a.e. $\theta\in X.$ 
\end{lemma}
\begin{proof}
Also routine.
\end{proof}

Our next unitary transformation is a simple unitary operator on the Hilbert space 
$L^2(\Omega\times\Pi\times\mathcal{L}^*)$.
Recall that, for each $L\in\Pi$, $x_L\in\R^n$ was selected so that $[x_L,L]\in\Gamma$. Define a function
$w_2:\Omega\times\Pi\times\mathcal{L}^*\to\T$ by $w_2(\omega,L,\nu)=e^{-2\pi i\nu\cdot x_L}$, for
all $(\omega,L,\nu)\in\Omega\times\Pi\times\mathcal{L}^*$. Pointwise multiplication by such a $\T$-valued
continuous function is a unitary operator and we formulate this in a lemma.

\begin{lemma}\label{W2}
For each $g\in L^2(\Omega\times\Pi\times\mathcal{L}^*)$, let $W_2g=w_2g$. That is,
\[
(W_2g)(\omega,L,\nu)=
e^{-2\pi i\nu\cdot x_L}g(\omega,L,\nu),
\]
 for a.e. $(\omega,L,\nu)\in\Omega\times\Pi\times\mathcal{L}^*$. Then $W_2$ is a unitary operator on 
$L^2(\Omega\times\Pi\times\mathcal{L}^*)$ and $W_2^{-1}g=\overline{w_2}g$, for all 
$g\in L^2(\Omega\times\Pi\times\mathcal{L}^*)$.
\end{lemma}

\begin{rem}
The function $w_2$ plays an important role in the characterization of $\Gamma$--shift-invariant closed
subspaces of $L^2(\R^n)$ obtained in Section 6, 
so it is important to note that this function is independent of the choice of the 
vector system associated with $\Gamma$. Indeed, 
if $\{x'_L:L\in\Pi\}$ is another vector system associated with $\Gamma$, then, for each $L\in\Pi$,
 $[x_L,L]^{-1}[x'_L,L]=[x'_L-x_L,{\rm id}]\in T$, so
$k=x'_L-x_L\in\mathcal{L}$.
Thus, $e^{-2\pi i\nu\cdot x'_L}=e^{-2\pi i\nu\cdot (x_L+k)}=e^{-2\pi i\nu\cdot x_L}$, for all $\nu\in\mathcal{L}^*$.
\end{rem}

Let $\{\delta_\nu:\nu\in\mathcal{L}^*\}$ be the orthonormal basis of $\ell^2(\mathcal{L}^*)$ consisting
of the usual delta functions.  We will use $L^2\big(\Pi\Omega,\ell^2(\mathcal{L}^*)\big)$, which is a
separable Hilbert space.
As we work with each of $L^2(\R^n)$, $L^2(\Omega\times\Pi\times\mathcal{L}^*)$, 
$\ell^2(\mathcal{L}^*)$, and
$L^2\big(\Pi\Omega,\ell^2(\mathcal{L}^*)\big)$, any inner product that arises will be
denoted $\langle\cdot,\cdot\rangle$ and the corresponding norm denoted by $\|\cdot\|$,
relying on the reader to know which Hilbert space this inner product or norm belongs
to from the context.

We now turn to the definition of a unitary transformation of $L^2(\Omega\times\Pi\times\mathcal{L}^*)$
into $L^2\big(\Pi\Omega,\ell^2(\mathcal{L}^*)\big)$. As above, for any $\theta\in\Pi\Omega$, there are
unique $\omega_\theta\in\Omega$ and $L_\theta\in\Pi$ so that $L_\theta\omega_\theta=\theta$. For any
$g\in L^2(\Omega\times\Pi\times\mathcal{L}^*)$ and $\nu\in\mathcal{L}^*$, let 
$g_\nu(\theta)=g(\omega_\theta,L_\theta,\nu)$, for a.e. $\theta\in\Pi\Omega$. Then 
$g_\nu\in L^2(\Pi\Omega)$, for each $\nu\in\mathcal{L}^*$, and 
$\|g\|^2=\sum_{\nu\in\mathcal{L}^*}\|g_\nu\|^2$. Define $W_3g:\Pi\Omega\to\ell^2(\mathcal{L}^*)$ by
\begin{equation}\label{W3def}
W_3g(\theta)=\sum_{\nu\in\mathcal{L}^*}g_\nu(\theta)\delta_\nu, \text{ for a.e. }\theta\in\Pi\Omega.
\end{equation}
\begin{lemma}\label{W3}
If $W_3g$ is defined by \eqref{W3def}, then $W_3g\in L^2\big(\Pi\Omega,\ell^2(\mathcal{L}^*)\big)$, 
for each $g\in L^2(\Omega\times\Pi\times\mathcal{L}^*)$ and $W_3$ is a unitary map of 
$L^2(\Omega\times\Pi\times\mathcal{L}^*)$
onto $L^2\big(\Pi\Omega,\ell^2(\mathcal{L}^*)\big)$. Moreover, for $F\in 
L^2\big(\Pi\Omega,\ell^2(\mathcal{L}^*)\big)$, 
\begin{equation}\label{W31}
W_3^{-1}F(\omega,L,\nu)=\langle F(L\omega),\delta_\nu\rangle,
\end{equation}
for a.e. 
$(\omega,L,\nu)\in \Omega\times\Pi\times\mathcal{L}^*$.
\end{lemma}
\begin{proof}
For each $\kappa\in\mathcal{L}^*$, the map 
$\theta\mapsto\langle W_3g(\theta),\delta_\kappa\rangle=g_\kappa(\theta)$ is measurable. Since 
$\{\delta_\kappa:\kappa\in\mathcal{L}^*\}$ is a countable basis of $\ell^2(\mathcal{L}^*)$, $W_3g$ is
measurable. Moreover, using \eqref{W3def}
\begin{equation}\label{W3sqsum}
\int_{\Pi\Omega}\|W_3g(\theta)\|^2d\theta=\int_{\Pi\Omega}\sum_{\nu\in\mathcal{L}^*}|g_\nu(\theta)|^2d\theta
=\sum_{\nu\in\mathcal{L}^*}\int_{\Pi\Omega}|g_\nu(\theta)|^2d\theta=\|g\|^2.
\end{equation}
Therefore, $W_3g\in L^2\big(\Pi\Omega,\ell^2(\mathcal{L}^*)\big)$ and $W_3$ is clearly linear. By
\eqref{W3sqsum}, $W_3$ is an isometry. If $F\in L^2\big(\Pi\Omega,\ell^2(\mathcal{L}^*)\big)$ and 
\eqref{W31} is used to define $W_3^{-1}F$ on $\Omega\times\Pi\times\mathcal{L}^*$, then one verifies,
with an argument similar to \eqref{W3sqsum} in reverse, that $W_3^{-1}F$ is square-integrable on
$\Omega\times\Pi\times\mathcal{L}^*$ and $W_3W_3^{-1}F=F$.
\end{proof}

\section{The natural representation}
We continue with the notation of the previous section.
The action of $\Gamma$ on $\R^n$ generates a unitary representation. 
For $[x,M]\in\Gamma$, $\pi[x,M]$ is the unitary operator given by, for $f\in L^2(\R^n)$,
\[
\pi[x,M]f(z)=f([x,M]^{-1}z)=f(M^{-1}z-x), \text{ for a.e. }z\in\R^n.
\]
It is routine to show that $\pi$ is a homomorphism of $\Gamma$ into $\mathcal{U}\big(L^2(\R^n)\big)$.
As mentioned in the introduction, we will call $\pi$ the {\em natural representation} of $\Gamma$.
The Fourier transform, as a unitary map $\mathcal{F}$ on $L^2(\R^n)$, intertwines $\pi$ with an equivalent representation
$\widehat{\pi}$ by $\widehat{\pi}[\ell,M]=\mathcal{F}\pi[\ell,M]\mathcal{F}^{-1}$. A brief calculation, using
the fact that $M^{-1}$ agrees with the transpose matrix of $M\in{\rm O}_n$,
shows that, for $\xi\in L^2(\R^n)$,
\begin{equation}\label{pihat}
\widehat{\pi}[x,M]\xi(\omega)=e^{2\pi i(M^{-1}\omega)\cdot x}\xi(M^{-1}\omega), 
\text{ for a.e. } \omega\in\R^n.
\end{equation}

Recall from Lemma \ref{W1} that $W_1\xi(\omega,L,\nu)=\xi\big(L(\omega+\nu)\big)$, for almost every
$(\omega,L,\nu)\in 
\Omega\times\Pi\times\mathcal{L}^*$, for  each $\xi\in L^2(\R^n)$, defines a unitary map
of $L^2(\R^n)$ onto $L^2(\Omega\times\Pi\times\mathcal{L}^*)$.
Let $\widehat{\pi}_1$
be the representation of $\Gamma$ on $L^2(\Omega\times\Pi\times\mathcal{L}^*)$ given by
$\widehat{\pi}_1[x,M]=W_1\widehat{\pi}[x,M]W_1^{-1}$, for $[x,M]\in\Gamma$. For 
$g\in L^2(\Omega\times\Pi\times\mathcal{L}^*)$, let $\xi=W_1^{-1}g$. Then, for $[x,M]\in\Gamma$ and
a.e. $\omega\in \Omega$, $L\in\Pi$, and $\nu\in\mathcal{L}^*$.
\[
\widehat{\pi}_1[x,M]g(\omega,L,\nu)=\widehat{\pi}[x,M]\xi\big(L(\omega+\nu)\big)=
e^{2\pi i(M^{-1}L(\omega+\nu))\cdot x}\xi\big(M^{-1}L(\omega+\nu)\big).
\]
Setting $\theta=M^{-1}L(\omega+\nu)$, note that $\omega_\theta=\omega$, $L_\theta=M^{-1}L$, and
$\nu_\theta=\nu$. 
Thus, 
\[
\widehat{\pi}_1[x,M]g(\omega,L,\nu)=e^{2\pi i(M^{-1}L(\omega+\nu))\cdot x}
g(\omega,M^{-1}L,\nu), 
\]
for a.e. $(\omega,L,\nu)\in \Omega\times\Pi\times\mathcal{L}^*$. Since $[x,M]\in \Gamma$, there exists 
$\ell\in
\mathcal{L}$ so that $x=\ell+x_M$. For any $\nu\in\mathcal{L}^*$ and $N\in\Pi$, $e^{2\pi i(N\nu)\cdot\ell}=1$.
Therefore, $e^{2\pi i(M^{-1}L(\omega+\nu))\cdot x}$ can be written as 
$e^{2\pi i(M^{-1}L\nu)\cdot x_M}e^{2\pi i(M^{-1}L\omega)\cdot x}=
e^{2\pi i\nu\cdot (L^{-1}Mx_M)}e^{2\pi i(M^{-1}L\omega)\cdot x}$
and
\begin{equation}\label{pi1first}
\widehat{\pi}_1[x,M]g(\omega,L,\nu)=e^{2\pi i\nu\cdot (L^{-1}Mx_M)}e^{2\pi i(M^{-1}L\omega)\cdot x}
g(\omega,M^{-1}L,\nu). 
\end{equation}
We can use the rather elementary, but strategic, unitary $W_2$ from Lemma \ref{W2} to eliminate
the first factor in the right hand side of \eqref{pi1first}.
Recall that 
$W_2g=w_2g$, for all $g\in L^2(\Omega\times\Pi\times\mathcal{L}^*)$, where
$w_2(\omega,L,\nu)=e^{-2\pi i\nu\cdot x_L}$, 
for $(\omega,L,\nu)\in \Omega\times\Pi\times\mathcal{L}^*$.
Let $\widehat{\pi}_2[x,M]=W_2\widehat{\pi}_1[x,M]W_2^{-1}$, for all $[x,M]\in\Gamma$. Then, for
$g\in L^2(\Omega\times\Pi\times\mathcal{L}^*)$,
\begin{align*}
&\widehat{\pi}_2[x,M]g(\omega,L,\nu) = e^{-2\pi i\nu\cdot x_L}\widehat{\pi}_1[x,M]W_2^{-1}g(\omega,L,\nu) \\
&= e^{-2\pi i\nu\cdot x_L}e^{2\pi i\nu\cdot (L^{-1}Mx_M)}e^{2\pi i(M^{-1}L\omega)\cdot x}
W_2^{-1}g(\omega,M^{-1}L,\nu)\\
&= e^{-2\pi i\nu\cdot x_L}e^{2\pi i\nu\cdot (L^{-1}Mx_M)}e^{2\pi i(M^{-1}L\omega)\cdot x}
 e^{2\pi i\nu\cdot x_{M^{-1}L}}g(\omega,M^{-1}L,\nu).
\end{align*}
The scalar in the previous line can be rearranged to
\[
 e^{2\pi i\nu\cdot (L^{-1}Mx_M+x_{M^{-1}L}-x_L  )}e^{2\pi i(L\omega)\cdot (Mx)}.
\]
Consider the expression $L^{-1}Mx_M+x_{M^{-1}L}-x_L$. By Lemma \ref{cocycle} with $K=M$ and $N=M^{-1}L$,
we have that $L^{-1}Mx_M+x_{M^{-1}L}-x_L=\ell$, for some $\ell\in\mathcal{L}$. Since $\nu\in\mathcal{L}^*$,
$e^{2\pi i\nu\cdot (L^{-1}Mx_M+x_{M^{-1}L}-x_L  )}=1$. Finally, for $[x,M]\in\Gamma$ and
$g\in L^2(\Omega\times\Pi\times\mathcal{L}^*)$,
\begin{equation}\label{pi2}
\widehat{\pi}_2[x,M]g(\omega,L,\nu)=e^{2\pi i(L\omega)\cdot (Mx)}g(\omega,M^{-1}L,\nu)
\end{equation}
for a.e. $\omega\in\Omega$, $L\in\Pi$ and $\nu\in\mathcal{L}^*$. 

The next step is to conjugate by the unitary $W_3:L^2(\Omega\times\Pi\times\mathcal{L}^*)\to
L^2\big(\Pi\Omega,\ell^2(\mathcal{L}^*)\big)$ from Lemma \ref{W3}.
Let $\widehat{\pi}_3[x,M]=W_3\widehat{\pi}_2[x,M]W_3^{-1}$, for all $[x,M]\in\Gamma$. 
For $F\in L^2\big(\Pi\Omega,\ell^2(\mathcal{L}^*)\big)$,
\begin{eqnarray}\label{pi3_1}
\widehat{\pi}_3[x,M]F(\theta) & = &W_3\widehat{\pi}_2[x,M]W_3^{-1}F(\theta)\nonumber\\
& = &
\sum_{\nu\in\mathcal{L}^*}\left(\widehat{\pi}_2[x,M]W_3^{-1}F(\omega_\theta,L_\theta,\nu)\right)\delta_\nu
\nonumber\\
& = &
e^{2\pi i(L_\theta\omega_\theta)\cdot (Mx)}\sum_{\nu\in\mathcal{L}^*}
W_3^{-1}F(\omega_\theta,M^{-1}L_\theta,\nu)\delta_\nu.
\end{eqnarray}
Using \eqref{W31} we get $\sum_{\nu\in\mathcal{L}^*}
W_3^{-1}F(\omega_\theta,M^{-1}L_\theta,\nu)\delta_\nu=F(M^{-1}L_\theta\omega_\theta)$. Recall that
$L_\theta\omega_\theta=\theta$ and $M^t=M^{-1}$, so \eqref{pi3_1} implies
\begin{equation}\label{pi3}
\widehat{\pi}_3[x,M]F(\theta)=e^{2\pi i(M^{-1}\theta)\cdot x}F(M^{-1}\theta), \text{ for a.e. }\theta\in
\Pi\Omega,
\end{equation}
and each $[x,M]\in\Gamma$ and any $F\in L^2\big(\Pi\Omega,\ell^2(\mathcal{L}^*)\big)$. 
We summarize these considerations in a proposition.

\begin{prop}\label{pi2equiv}
If $\widehat{\pi}_2[x,M]$ is defined by \eqref{pi2} and $\widehat{\pi}_3[x,M]$ is defined by \eqref{pi3}, 
for each $[x,M]\in\Gamma$, then $\widehat{\pi}_2$ and $\widehat{\pi}_3$
are unitary representations of $\Gamma$ that are each  unitarily equivalent to the natural representation $\pi$.
\end{prop}

\section{Invariant closed subspaces}

\begin{defn}\label{L_inv_range}
Let $J:\Pi\Omega\to Gr\big(\ell^2(\mathcal{L}^*)\big)$ be a measurable range function. 
We say $J$ is $\Pi$-invariant, if $J(L\theta)=
J(\theta)$, for a.e. $\theta\in\Pi\Omega$ and every $L\in\Pi$.
\end{defn}

\begin{prop}\label{inv_pi3}
Let $\mathcal{K}$ be a closed subspace of $L^2\big(\Pi\Omega,\ell^2(\mathcal{L}^*)\big)$. Then 
$\mathcal{K}$ is $\widehat{\pi}_3$-invariant if and only if there exists a 
$\Pi$-invariant measurable range function, $J$, for
$\ell^2(\mathcal{L}^*)$ based on $\Pi\Omega$ such that $\mathcal{K}=\mathcal{M}_J$.
\end{prop}
\begin{proof}
First, suppose $J$ is a $\Pi$-invariant measurable range function for
$\ell^2(\mathcal{L}^*)$ based on $\Pi\Omega$ and $\mathcal{K}=\mathcal{M}_J$.
 Let $[x,M]\in\Gamma$. For any $F\in\mathcal{M}_J$, we know that
$F(\theta)\in J(\theta)$, for a.e. $\theta\in\Pi\Omega$, and $J$ is $\Pi$-invariant. This implies
\[
\widehat{\pi}_3[x,M]F(\theta)=e^{2\pi i\theta\cdot(Mx)}F(M^{-1}\theta)\in J(M^{-1}\theta)=J(\theta),
\]
for a.e. $\theta\in\Pi\Omega$. Thus, $\widehat{\pi}_3[x,M]F\in \mathcal{K}$, for all $F\in\mathcal{K}$ and 
$[x,M]\in\Gamma$. That is, $\mathcal{K}$ is $\widehat{\pi}_3$-invariant.

Now suppose $\mathcal{K}$ is a $\widehat{\pi}_3$-invariant closed subspace of
$L^2\big(\Pi\Omega,\ell^2(\mathcal{L}^*)\big)$. By Remark \ref{determing2},
$\mathcal{D}=\{e_\ell|_{_{\Pi\Omega}}:\ell\in\mathcal{L}\}$ is a determining set for $L^1(\Pi\Omega)$, where 
$e_\ell(\theta)=e^{2\pi i\theta\cdot\ell}$, for all $\theta\in \R^n$.
If $F\in \mathcal{K}$ and $\ell\in\mathcal{L}$, then $[\ell,{\rm id}]\in T\subseteq\Gamma$ and \eqref{pi3}
says
\[
\widehat{\pi}_3[\ell,{\rm id}]F(\theta)=e^{2\pi i\theta\cdot\ell}F(\theta)=e_\ell(\theta)F(\theta), 
\text{ for a.e. }\theta\in\Pi\Omega.
\]
Therefore, $e_\ell F\in\mathcal{K}$, for each $\ell\in\mathcal{L}$. In the language of 
\cite{BownikRoss}, $\mathcal{K}$ is multiplicatively-invariant with respect to $\mathcal{D}$, denoted
$\mathcal{D}-MI$. By Theorem 2.4 of \cite{BownikRoss}, there exists a measurable range function
 $J:\Pi\Omega\to Gr\big(\ell^2(\mathcal{L}^*)\big)$ so that 
\[
\mathcal{K}=\mathcal{M}_{J}=
\{F\in L^2\big(\Pi\Omega,\ell^2(\mathcal{L}^*)\big):F(\theta)\in J(\theta), \text{ for a.e. }\theta\in\Pi\Omega\}.
\]
Moreover, $J$ is unique up to almost everywhere agreement and, if $\mathcal{A}$ is any countable dense subset 
of $\mathcal{K}$, $J(\theta)$ is the closed linear span of $\{\varphi(\theta):\varphi\in\mathcal{A}\}$, for
a.e. $\theta\in\Pi\Omega$. Fix such an $\mathcal{A}$.

Fix $L\in\Pi$ so that $[x_L,L]\in\Gamma$ and use
\eqref{pi3} again to get
\begin{equation}\label{pi3xL}
\widehat{\pi}_3[x_L,L]F(L\omega)=e^{2\pi i\omega\cdot(x_L)}F(\omega), \text{ for a.e. }\omega\in
\Omega.
\end{equation}
%Note, for $v\in\ell^2(\mathcal{L}^*)$, $v\in J(\omega)$ if and only if $e^{2\pi i\omega\cdot(x_L)}v\in J(\omega)$.
For $\varphi\in\mathcal{A}$, let $\varphi^L=\widehat{\pi}_3[x_L,L]\varphi$ and let
$L\cdot\mathcal{A}=\{\varphi^L:\varphi\in\mathcal{A}\}$. Then $L\cdot\mathcal{A}$ is a dense subset
of $\mathcal{K}$. Therefore, by \eqref{pi3xL}, for a.e. $\omega\in\Omega$, 
\[
J(L\omega)=\overline{\rm span}\{\varphi^L(L\omega):\varphi\in\mathcal{A}\}=
\overline{\rm span}\{\varphi(\omega):\varphi\in\mathcal{A}\}=J(\omega).
\]
Since $L\in\Pi$ is arbitrary, $J$ is $\Pi$-invariant.
\end{proof}
Before formulating the main result, we recall the various ingredients: $\Gamma$ is a crystal group in dimension
$n$ with point group $\Pi$ and associated lattice $\mathcal{L}$. The dual lattice of $\mathcal{L}$ is
$\mathcal{L}^*$. There is an open subset $\Omega$ of $\R^n$ described in Proposition \ref{Rproperties}
and the collection of sets $\{L(\Omega+\nu):L\in\Pi,\nu\in\mathcal{L}^*\}$ tiles $\R^n$ up to a null set.
\begin{theorem}\label{main_theorem}
Let $\mathcal{V}$ be a closed subspace of 
$L^2(\R^n)$. Then $\mathcal{V}$ is invariant under shifts by elements of $\Gamma$ if and only if there exists 
a $\Pi$-invariant measurable range function $J$ for $\ell^2(\mathcal{L}^*)$ based on $\Pi\Omega$ such that,
for any $f\in\mathcal{V}$, there exists $F\in\mathcal{M}_J$ so that, for almost every
$\omega\in\Omega$, every $L\in\Pi$,
and every $\nu\in\mathcal{L}^*$,
\begin{equation}\label{main}
\widehat{f}\big(L(\omega+\nu)\big)=e^{2\pi i\nu\cdot x_L}\langle F(L\omega),\delta_\nu\rangle.
\end{equation}
\end{theorem}
\begin{proof}
Suppose that $\mathcal{V}$ is a closed subspace of $L^2(\R^n)$ that is 
invariant under shifts by elements of $\Gamma$. That is, $\mathcal{V}$ is $\pi$-invariant.
Let $U=W_3\circ W_2\circ W_1\circ\mathcal{F}$, a unitary map that intertwines the natural 
representation $\pi$ with $\widehat{\pi}_3$. Therefore, $U\mathcal{V}$ is a $\widehat{\pi}_3$-invariant
closed subspace of $L^2\big(\Pi\Omega,\ell^2(\mathcal{L}^*)\big)$. By Proposition \ref{inv_pi3}, there
exists a $\Pi$-invariant measurable range function $J:\Pi\Omega\to Gr\big(\ell^2(\mathcal{L}^*)\big)$
such that $U\mathcal{V}=\mathcal{M}_J$. For any $f\in\mathcal{V}$, let $F=Uf$. This can be written as
$W_1\widehat{f}=W_2^{-1}W_3^{-1}F$, which is an element of $L^2(\Omega\times\Pi\times\mathcal{L}^*)$.
By the definition of $W_1$ given just before Lemma \ref{W1}, 
\begin{equation}\label{LHS}
W_1\widehat{f}(\omega,L,\nu)=\widehat{f}\big(L(\omega+\nu)\big), 
\end{equation}
for a.e. $(\omega,L,\nu)\in\Omega\times\Pi\times\mathcal{L}^*$.
On the other hand, using the equations for $W_2^{-1}$ and $W_3^{-1}$ given in Lemmas \ref{W2} and \ref{W3},
we have
\begin{equation}\label{RHS}
W_2^{-1}W_3^{-1}F(\omega,L,\nu)=e^{2\pi i\nu\cdot x_L}\langle F(L\omega),\delta_\nu\rangle,
\end{equation}
for a.e. $(\omega,L,\nu)\in\Omega\times\Pi\times\mathcal{L}^*$. Comparing \eqref{LHS} with \eqref{RHS}
yields \eqref{main}.

Conversely, suppose $J$ is a $\Pi$-invariant measurable range function for $\ell^2(\mathcal{L}^*)$ 
based on $\Pi\Omega$ and that $\mathcal{V}$ is related to $\mathcal{M}_J$ via \eqref{main}.
Then $\mathcal{M}_J$ is $\widehat{\pi}_3$-invariant by Proposition \ref{inv_pi3} and unpacking \eqref{main}
shows that $U^{-1}\mathcal{M}_J=\mathcal{V}$. Since $U$ is the intertwining unitary between $\pi$ and
$\widehat{\pi}_3$, $\mathcal{V}$ must be $\pi$-invariant. Thus, $\mathcal{V}$ is invariant under shifts 
by elements of $\Gamma$.
\end{proof}
Theorem \ref{main_theorem} is a direct generalization of Proposition 1.5 in \cite{Bownik}, where it is attributed to
Helson \cite{Helson}. It is not surprising that the range function associated to a $\pi$-invariant closed
subspace must be $\Pi$-invariant. However, the appearance of the factor $e^{2\pi i\nu\cdot x_L}$ in
\eqref{main} is not so obvious. As noted earlier, for a given $L\in\Pi$, the choice of $x_L$ such that
$[x_L,L]\in\Gamma$ is not unique. However, if $y\in\R^n$ also satisfies $[y,L]\in\Gamma$, then 
$x_L-y\in\mathcal{L}$. Thus, $e^{2\pi i\nu\cdot y}=e^{2\pi i\nu\cdot x_L}$, for all $\nu\in\mathcal{L}^*$.
To help understand the role this term plays, we present an example of
a  group, necessarily non-symmorphic, where this factor is nontrivial.

\begin{ex}\label{pgex}
The patch of brick wall illustrated in Figure 1 is meant to be a region of a pattern extending in all directions. 
The symmetry group of this pattern is often
denoted $pg$, so set $\Gamma=pg$. 
If the origin is placed at the bottom of one of the vertical line segments, let ${\bf u}$ denote a horizontal
vector pointing right whose length equals the length of one brick. Let ${\bf v}$ be an upward pointing vector whose length equals twice the width of a brick. Then, we can take the translation lattice of $\Gamma$ to
be $\mathcal{L}=\{j{\bf u}+k{\bf v}: j,k\in\Z\}$. Let $\sigma\in {\rm O}_2$ denote the reflection about the
horizontal axis. This is not a symmetry of the pattern, but $[\frac{1}{2}{\bf u},\sigma]\in \Gamma$. Indeed,
the point group of $\Gamma$ is $\Pi=\{{\rm id},\sigma\}$, isomorphic to $\Z_2$. 
If $T=\{[j{\bf u}+k{\bf v},{\rm id}]: j,k\in\Z\}$
is the translation subgroup of $\Gamma$, then $\{{\rm id} \mapsto 0;\  \sigma \mapsto \frac{1}{2}{\bf u}\}$ is a vector system associated with $\Gamma$.
		\begin{figure}[H]
			\centering
			\begin{tikzpicture}[scale=.6]
				\foreach \i [evaluate=\i as \y using \i]  in {0,2}
				\foreach \x in {-4,-2,...,4}
				{
					\shade[thick, lower right = white!90!black, upper left = white!40!black] (\x, \y) -- (\x+2, \y) -- (\x+2,\y+1) -- (\x, \y+1);		
					\draw (\x, \y) -- (\x+2, \y) -- (\x+2,\y+1) -- (\x, \y+1) -- (\x,\y) -- (\x+2,\y);
					\draw (\x+.8, \y+.3) -- (\x+.8, \y+.7) -- (\x+1.2,\y+.7); 
				}
				\foreach \i [evaluate=\i as \y using \i]  in {1,3}
				\foreach \x in {-4,-2,...,4}
				{
					\shade[thick, lower left = white!40!black, upper right = white!90!black]  (\x+1, \y) -- (\x+3, \y) -- (\x+3,\y+1) -- (\x+1, \y+1);	
					\draw (\x+1, \y) -- (\x+3, \y) -- (\x+3,\y+1) -- (\x+1, \y+1) -- (\x+1,\y) -- (\x+3,\y);	
					\draw (\x+1.8,\y+.7) -- (\x+1.8, \y+.3) -- (\x+2.2, \y+.3); 
				}
\draw[dotted,thick] (-5.5,2)--(-4.5,2);	
\draw[dotted,thick] (7.5,2)--(8.5,2);	
\draw[dotted, thick] (1,4.25)--(1,4.75);	
\draw[dotted, thick] (1,-.25)--(1,-.75);		
			\end{tikzpicture}
			\caption{A pattern illustrating the symmetries of the wallpaper group $pg$.}
		\end{figure}
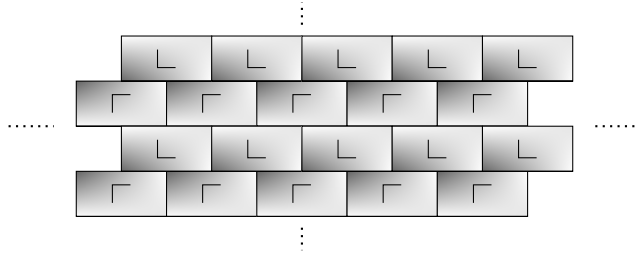
Fix the basis for $\R^2$ to be $\{{\bf u}, {\bf v}\}$ and represent elements $x\in\R^2$ as $x=\begin{pmatrix}
x_1\\
x_2
\end{pmatrix}$. Then $\mathcal{L}=\left\{\ell=\begin{pmatrix}
\ell_1\\
\ell_2
\end{pmatrix}:\ell_1,\ell_2\in\Z\right\}=\Z^2$, viewed as a lattice in $\R^2$. So $\mathcal{L}^*=\Z^2$, as well.
We set $\Gamma^*=\mathcal{L}^*\rtimes\Pi$.
Let $\omega_0=\begin{pmatrix}
0\\
1/4
\end{pmatrix}$ and form $\Omega$ as in Definition \ref{Dirichlet}. Then 
$\Omega=\left\{\omega=\begin{pmatrix}
\omega_1\\
\omega_2
\end{pmatrix};-\frac{1}{2}<\omega_1<\frac{1}{2}, 0<\omega_2<\frac{1}{2}\right\}$ as illustrated in
Figure 2. The open set $\sigma\Omega$ is outlined as well. Note that $\Pi\Omega=\Omega\cup \sigma\Omega$.

\begin{figure}[H]
			\centering
			\begin{tikzpicture}[scale=1.8]
\draw[->] (-2.5,0)--(2.5,0);
\draw[->] (0,-.85)--(0,1);
\path[pattern=north west lines] (-.5,0) rectangle (.5,.5);
\draw[fill] (0,.25) circle [radius=0.02];
\draw[fill] (1,.25) circle [radius=0.02];
\draw[fill] (2,.25) circle [radius=0.02];
\draw[fill] (-1,.25) circle [radius=0.02];
\draw[fill] (-2,.25) circle [radius=0.02];
\draw[fill] (0,.75) circle [radius=0.02];
\draw[fill] (1,.75) circle [radius=0.02];
\draw[fill] (2,.75) circle [radius=0.02];
\draw[fill] (-1,.75) circle [radius=0.02];
\draw[fill] (-2,.75) circle [radius=0.02];
\draw[fill] (0,-.25) circle [radius=0.02];
\draw[fill] (1,-.25) circle [radius=0.02];
\draw[fill] (2,-.25) circle [radius=0.02];
\draw[fill] (-1,-.25) circle [radius=0.02];
\draw[fill] (-2,-.25) circle [radius=0.02];
\draw[fill] (0,-.75) circle [radius=0.02];
\draw[fill] (1,-.75) circle [radius=0.02];
\draw[fill] (2,-.75) circle [radius=0.02];
\draw[fill] (-1,-.75) circle [radius=0.02];
\draw[fill] (-2,-.75) circle [radius=0.02];
\node at (.25,.35) {$\Omega$};
\draw[dashed] (-.5,-.5) rectangle (.5,0);
\node at (.25,-.35) {$\sigma\Omega$};
\end{tikzpicture}
			\caption{The $\Gamma^*$-orbit of $\omega_0$ and domain $\Omega$.}
		\end{figure}
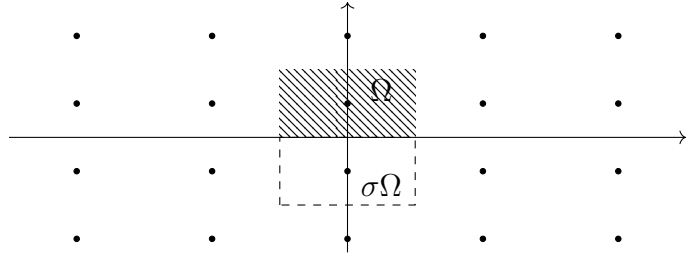
Let's construct a simple $\Pi$-invariant measurable range function for $\ell^2(\mathcal{L}^*)$ based on
$\Pi\Omega$. Let $\kappa=\begin{pmatrix}
1\\
0
\end{pmatrix}\in\mathcal{L}^*$ and consider the one-dimensional 
subspace $V=\C(\delta_0+\delta_\kappa)=\{\alpha(\delta_0+\delta_\kappa):\alpha\in\C\}$ of $\ell^2(\mathcal{L}^*)$.
Let $E$ be a measurable subset of $\Omega$ as pictured in Figure 3 where we have suppressed the
vertical axis and changed the scale a little.

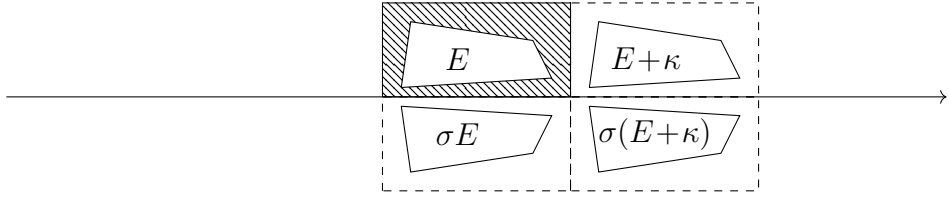
\begin{figure}[H]
			\centering
			\begin{tikzpicture}[scale=2.5]
\draw[->] (-2.5,0)--(2.5,0);
\draw[pattern=north west lines] (-.5,0) rectangle (.5,.5);
\draw[dashed] (-.5,-.5) rectangle (.5,0);
\draw[dashed] (.5,0) rectangle (1.5,.5);
\draw[dashed] (.5,-.5) rectangle (1.5,0);
\draw[fill=white] (-.4,.05)--(.4,.1)--(.3,.3)--(-.35,.4)--(-.4,.05);
\draw (.6,.05)--(1.4,.1)--(1.3,.3)--(.65,.4)--(.6,.05);
\draw (-.4,-.05)--(.4,-.1)--(.3,-.3)--(-.35,-.4)--(-.4,-.05);
\draw (.6,-.05)--(1.4,-.1)--(1.3,-.3)--(.65,-.4)--(.6,-.05);
\node at (-.1,.2) {$E$};
\node at (-.1,-.2) {$\sigma E$};
\node at (.9,.2) {$E\!+\!\kappa$};
\node at (.95,-.2) {$\sigma(E\!+\!\kappa)$};
\end{tikzpicture}
			\caption{The measurable subset $E$ of $\Omega$ and three of its shifts.}
		\end{figure}
Now consider the map $J:\Pi\Omega\to Gr\big(\ell^2(\mathcal{L}^*)\big)$ given by $J(\theta)=V$, for
$\theta\in E\cup \sigma E$, and $J(\theta)=\{0\}$, for $\theta\in\Pi\Omega\setminus(E\cup \sigma E)$. Then,
$\mathcal{M}_J$ consists of all $F\in L^2\big(\Pi\Omega,\ell^2(\mathcal{L}^*)\big)$ such that $F(\theta)=0$,
for a.e. $\theta\in\Pi\Omega\setminus(E\cup \sigma E)$, and $F(\theta)\in\C(\delta_0+\delta_\kappa)$, for
a.e. $\theta\in E\cup \sigma E$.
Before looking at 
\eqref{main} in this situation, note that $x_{\rm id}=\begin{pmatrix}
0\\
0
\end{pmatrix}$ and $x_\sigma=\begin{pmatrix}
1/2\\
0
\end{pmatrix}$. So, $e^{2\pi i\nu\cdot x_{\rm id}}=1$ and $e^{2\pi i\nu\cdot x_\sigma}=(-1)^{\nu_1}$, for all
$\nu=\begin{pmatrix}
\nu_1\\
\nu_2
\end{pmatrix}\in\mathcal{L}^*$. Consider the following three properties for an $f\in L^2(\R^2)$:\\
\indent (a) $\widehat{f}(\omega+\kappa)=\widehat{f}(\omega)$, for a.e. $\omega\in\Omega$.\\
\indent (b) $\widehat{f}\big(\sigma(\omega+\kappa)\big)=
-\widehat{f}\big(\sigma(\omega)\big)$, for a.e. $\omega\in\Omega$.\\
\indent (c) $\widehat{f}(\theta)=0$, for all $\theta\in \R^2\setminus\big(
E\cup (E+\kappa)\cup \sigma E\cup \sigma(E+\kappa)\big)$.\\
What Theorem \ref{main_theorem} implies for the simple range function $J$ is that, if $\mathcal{V}$ denotes the set of all
$f\in L^2(\R^2)$ whose Fourier transform $\widehat{f}$ satisfies (a), (b), and (c), then $\mathcal{V}$ is
a closed subspace of $L^2(\R^2)$ that is invariant under shifts from the wallpaper group $pg$.
\end{ex}

Although the range function used in Example \ref{pgex} is simple compared to the complexity that is possible,
the example does suggest that the statement of Theorem \ref{main_theorem} can be refined. 
We return to 
$\Gamma$ denoting a crystal group in dimension $n$ with $\Pi$, $\mathcal{L}$, and $\Omega$ as before.
For each $L\in\Pi$, define a unitary operator $U_L$ on $\ell^2(\mathcal{L}^*)$ by
$U_Lh(\nu)=e^{2\pi i\nu\cdot x_L}h(\nu)$, for all $\nu\in\mathcal{L}^*$ and for each $h\in\ell^2(\mathcal{L}^*)$.
In Theorem \ref{main_theorem}, we used range functions based on $\Pi\Omega$ that are $\Pi$-invariant, so
completely determined by their restriction to $\Omega$. We now give a version where $\Pi$-invariance is replaced with
twisting by the unitaries $U_L$, $L\in\Pi$.
\begin{defn}
For a measurable range function $J:\Omega\to Gr\big(\ell^2(\mathcal{L}^*)\big)$, let $J^\Gamma:\Pi\Omega
\to Gr\big(\ell^2(\mathcal{L}^*)\big)$ be given by
$J^\Gamma(L\omega)=U_LJ(\omega)$, for each $\omega\in\Omega$ and all $L\in\Pi$.
\end{defn}
If $J$ is a measurable range function for $\ell^2(\mathcal{L}^*)$ based on $\Omega$, then it is clear that
$J^\Gamma$ is a measurable range function for $\ell^2(\mathcal{L}^*)$ based on $\Pi\Omega$. Note that we can
create a $\Pi$-invariant range function $J'$  from $J$ based on $\Pi\Omega$ by simply letting
$J'(L\omega)=J(\omega)$, for all $L\in\Pi$, $\omega\in\Omega$. We call $J'$ the $\Pi$-invariant extension 
of $J$. The following lemma is immediate from the definitions.
\begin{lemma}\label{JPi}
Let $J$ be a measurable range function for $\ell^2(\mathcal{L}^*)$ based on $\Omega$ and let $J'$ be
the $\Pi$-invariant extension. For any $G\in L^2\big(\Pi\Omega,\ell^2(\mathcal{L}^*)\big)$,  we have
$G(L\omega)\in J'(L\omega)$ if and only if $U_L\big(G(L\omega)\big)\in J^\Gamma(L\omega)$, for all $\omega\in\Omega$ and $L\in\Pi$.
\end{lemma}
We can now present a modified statement of Theorem \ref{main_theorem} as a corollary.
\begin{cor}\label{main_cor}
Let $\mathcal{V}$ be a closed subspace of $L^2(\R^n)$. Then $\mathcal{V}$ is invariant under shifts by 
elements of $\Gamma$ if and only if there exists a measurable range function $J$ for $\ell^2(\mathcal{L}^*)$
based on $\Omega$ such that, for any $f\in\mathcal{V}$, there exists $F\in \mathcal{M}_{J^\Gamma}$ so that
$\widehat{f}\big(L(\omega+\nu)\big)=\langle F(L\omega),\delta_\nu\rangle$,
for almost every $\omega\in\Omega$, every $L\in\Pi$, and every $\nu\in\mathcal{L}^*$.
\end{cor}
\begin{proof}
Suppose that $\mathcal{V}$ is invariant under shifts by $\Gamma$, and let $K$ be a $\Pi$-invariant range function for $\ell^2(\mathcal{L}^*)$ based on $\Pi\Omega$, as described in Theorem \ref{main_theorem}. Thus for any $f \in \mathcal{V}$, there exists $G\in \mathcal{M}_K$ such that $\widehat{f}\big(L(\omega+\nu)\big)=e^{2\pi i\nu\cdot x_L}\langle G(L\omega),\delta_\nu\rangle$ for a.e. $\omega \in \Omega$ and all $L\in \Pi$ and $\nu\in\mathcal{L}^*$. Let $J$ denote the restriction of $K$ to $\Omega$. Then $J'=K$. Define
$F\in L^2\big(\Pi\Omega,\ell^2(\mathcal{L}^*)\big)$ by $F(L\omega)=U_L\big(G(L\omega)\big)$, for a.e. 
$\omega\in\Omega$ and any  $L\in\Pi$. By Lemma \ref{JPi}, $F\in\mathcal{M}_{J^\Gamma}$. And
\begin{equation}\label{MainCorEqn}
\widehat{f}\big(L(\omega+\nu)\big)=e^{2\pi i\nu\cdot x_L}\langle G(L\omega),\delta_\nu\rangle=\langle U_L\big(G(L\omega)\big),\delta_\nu\rangle
=\langle F(L\omega),\delta_\nu\rangle,
\end{equation}
for a.e. $\omega\in\Omega$, any $L\in\Pi$, and any $\nu\in\mathcal{L}^*$.

Conversely, suppose that for any $f\in\mathcal{V}$, there exists $F\in \mathcal{M}_{J^\Gamma}$ such that
\[
\widehat{f}\big(L(\omega+\nu)\big) = \langle F(L\omega),\delta_\nu \rangle
\]
for a.e. $\omega$ and all $L\in\Pi$, $\nu\in\mathcal{L}^*$. Defining $G \in L^2\big(\Pi\Omega, \ell^2(\mathcal{L}^*)\big)$ by $G(L\omega) = U_L^{-1} F(L\omega)$ for a.e. $\omega\in\Omega$, it is immediate that $\langle F(L\omega),\delta_\nu\rangle = e^{2\pi i\nu\cdot x_L}\langle G(L\omega),\delta_\nu\rangle$, so that
\[
\widehat{f}\big(L(\omega+\nu)\big) = e^{2\pi i\nu\cdot x_L} \langle G(L\omega),\delta_\nu\rangle. 
\]
We show that $G \in \mathcal{M}_K$, where $K$ is a $\Pi$-invariant measurable range function for $\ell^2(\mathcal{L}^*)$ based on $\Pi \Omega$. Since $F\in \mathcal{M}_{J^\Gamma}$ by assumption, we have
\[
U_L G(L\omega) = F(L\omega) \in J^\Gamma(L\omega), 
\]
which by Lemma 6.6 is equivalent to $G(L\omega) \in J^\prime (L\omega)$ for a.e. $\omega$, where $J^\prime$ is the $\Pi$-invariant extension of $J$. Thus $G\in \mathcal{M}_{J^\prime}$. It now follows from Theorem 6.3 that 
$\mathcal{V}$ is invariant under shifts by $\Gamma$.
\end{proof}
%\bibliographystyle{plain}
%\bibliography{Draft1Bib}{}

\end{document}